\documentclass[runningheads]{llncs}
\usepackage{amsfonts}
\usepackage{amsmath}

\usepackage{hyperref}

\usepackage{pgf, tikz}
\usepackage{color}

\def\RR{{\mathbb R}}
\def\QQ{{\mathbb Q}}
\def\AA{{\mathbb A}}
\def\ZZ{{\mathbb Z}}
\def\strut{\vphantom{\Large(} }

\def\antic{\textit{antic}}
\def\arb{\textit{arb}}
\def\eantic{\textit{e-antic}}
\def\Flint{\textit{Flint}}

\begin{document}
\title{Algebraic polytopes in Normaliz}
%
%
\author{Winfried Bruns\orcidID{0000-0002-7081-2261}}
\authorrunning{W. Bruns}
%
\institute{Institut f\"ur Mathematik, Universit\"at Osnabr\"uck, 49069 Osnabr\"uck, Germany
\email{wbrunas@uos.de}\\
\url{http://www.home.uni-osnabrueck.de/wbruns/} }
\maketitle              
\begin{abstract}
We describe the implementation  of algebraic polyhedra in Normaliz. In addition to convex hull computation/vertex enumeration, Normaliz computes triangulations, volumes, lattice points, face lattices and automorphism groups. The arithmetic is based on the package \eantic\ by V.~Delecroix. 

\end{abstract}

Algebraic polytopes lacking a rational realization are among the first geometric objects encountered in high school geometry: at least one vertex  of an equilateral triangle in the plane has non-rational coordinates. Three of the five Platonic solids, namely the tetrahedron, the icosahedron and the dodecahedron are non-rational, and,  among the $4$-dimensional regular polytopes, the $120$-cell and the $600$-cell live outside the rational world. 

But algebraic polytopes do not only appear in connection with Coxeter groups. Other contexts include enumerative combinatorics \cite{rote}, Dirichlet domains of hyperbolic group actions \cite{DelPage}, $\operatorname{SL}(2,\RR)$-orbit closures in the moduli space of translation surfaces, and parameter spaces and perturbation polyhedra of cut-generating functions in integer programming.

\section{Real embedded algebraic number fields}

The notion of convexity is defined over any ordered field, not only over the rationals $\QQ$ or the reals $\RR$. \emph{Real embedded algebraic number fields} are subfields of the real numbers (and therefore ordered) that have finite dimension as a $\QQ$-vector space. It is well known that such a field $\AA$ has a primitive element, i.e., an element $a$ such that no proper subfield of $\AA$ contains $a$. The minimal polynomial of $a$ is the least degree  monic polynomial $\mu$ with coefficients in $\QQ$ such that $\mu(a)=0$. It is an irreducible polynomial, and $\dim_\QQ\AA=\deg \mu$. In particular, every element $b$ of $\AA$ has a unique  representation  $b= \alpha_{n-1}a^{n-1}+\dots+\alpha_1a +\alpha_0$ with $\alpha_{n-1},\dots,\alpha_0\in\QQ$, $n=\deg \mu$. The arithmetic in $\AA$ is completely determined by $\mu$: addition is the addition of polynomials and multiplication is that of polynomials followed by reduction modulo $\mu$. The multiplicative inverse can be computed by the extended Euclidean algorithm. The unique determination of the coefficients $\alpha_i$ allows one to decide whether $b=0$. Every element of $\AA$ can be written as the quotient of a polynomial expression $\alpha_{n-1}a^{n-1}+\dots+\alpha_1a +\alpha_0$ with $\alpha_i\in\ZZ$ for all $i$ and an integer denominator; this representation is used in the implementation.

However, the algebraic structure alone does not define an ordering of $\AA$. For example, $\sqrt 2$ and $\sqrt{-2}$ cannot be distinguished algebraically: there exists an automorphism of $\QQ[\sqrt 2]$ that exchanges them. For the ordering we must fix a real number $a$ whose minimal polynomial is $\mu$. (Note that not every algebraic number field has an embedding into $\RR$.) In order to decide whether $b>0$ for some $b\in\AA$ we need a floating point approximation to $b$ of controlled precision.

Normaliz \cite{Nmz} uses the package \eantic\ of V. Delecroix \cite{e-antic} for the arithmetic and ordering in real algebraic number fields. The algebraic operations are realized by functions taken from the package \antic\ of W. Hart and F. Johansson \cite{antic} (imported to \eantic) while the controlled floating point arithmetic is delivered by the package \arb\ of F. Johansson \cite{arb}. Both packages are based on W. Hart's \Flint\ \cite{Flint}.

In order to specify an algebraic number field, one chooses the minimal polynomial $\mu$ of $a$ and an interval $I$ in $\RR$ such that $\mu$ has a unique zero in $I$, namely $a$. An initial approximation to $a$ is computed at the start. Whenever the current precision of $b$ does not allow to decide whether $b>0$, first the approximation of $b$ is improved, and if the precision of $a$ is not sufficient, it is replaced by one with twice the number of correct digits.

\section{Polyhedra}

A subset $P\subset \RR^d$ is a \emph{polyhedron} if it is the intersection of finitely many affine halfspaces:
$$
P=\bigcap_{i=0}^s H_i^+,\qquad H_i^+=\{x:\lambda_i(x) \ge \beta_i  \}, \qquad i=1,\dots,s,
$$
where $\lambda_i$ is a linear form and $\beta_i\in\RR$. It is a \emph{cone} if one can choose $\beta_i=0$ for all $i$, and it is a \emph{polytope} if it is bounded.

By the theorem of Minkowski-Weyl-Motzkin \cite[1.C]{BG} one can equivalently describe polyhedra by ``generators'': there exist $c_1,\dots,c_t\in\RR^d$ and $v_1,\dots,v_u\in \RR^d$ such that
$$
P=C+Q
$$
where $C=\bigl\{c\gamma_1c_1+\dots+\gamma_tc_t:\gamma_i\in\RR,\gamma_i\ge 0 \bigr\}$ is the \emph{recession cone} and $Q=\bigl\{\kappa_1v_1+\dots+\kappa_uv_u: \kappa_i\in\RR,\kappa_i\ge 0,\sum \kappa_i=1  \bigr\}$ is a polytope. These two descriptions are often called \emph{H-representation} and \emph{V-representation}. The conversion from H to V is  \emph{vertex enumeration} and the opposite conversion is \emph{convex hull computation}.

For theoretical and computational reasons it is advisable to present a polyhedron $P$ as the intersection of a cone and a hyperplane. Let $C(P)$ be the \emph{cone over} $P$, i.e., the smallest cone containing $P\times\{1\} $, and $D=\{x: x_{d+1}=1\}$ the \emph{dehomogenizing hyperplane}. Then $P$ can be identified  with $C(P)\cap D$. After this step, convex hull computation and vertex enumeration are two sides of the same coin, namely the dualization of cones. 
	
In the definition of polyhedra and all statements following it, the field $\RR$ can be replaced by an arbitrary subfield (and even by an arbitrary ordered field), for example a real algebraic number field $\AA$. The smallest choice for $\AA$ is $\QQ$: for it we obtain the class of \emph{rational polyhedra}. For general $\AA$ we get \emph{algebraic polyhedra}.

For the terminology related to polyhedra and further details we refer the treader to \cite{BG}.

\section{Normaliz}
Normaliz tackles many computational problems for rational and algebraic polyhedra:
\begin{itemize}
\item dual cones: convex hulls and vertex enumeration
\item projections of cones and polyhedra
\item triangulations, disjoint decompositions and Stanley decompositions
\item  Hilbert bases of rational, not necessarily pointed cones
\item  normalizations of affine monoids (hence the name)
\item  lattice points of polytopes and (unbounded) polyhedra
\item  automorphisms (euclidean, integral, rational/algebraic, combinatorial)
\item  face lattices and f-vectors
\item  Euclidean and lattice normalized volumes of polytopes
\item Hilbert (or Ehrhart) series and (quasi) polynomials under $\ZZ$-gradings
\item generalized (or weighted) Ehrhart series and Lebesgue integrals of polynomials over rational polytopes
\end{itemize}

Of course, not all of these computation goals make sense for algebraic polyhedra. The main difference between the rational and the non-rational case can be described as follows: the monoid of lattice points in a full dimensional cone is finitely generated if and only if the cone is rational.

Normaliz is based on a templated C++ library. The template allows one to choose the arithmetic, and so it would be possible to extend Normaliz to more general ordered fields. The main condition is that the arithmetic of the field has been coded in a C++ class library.  There is no restriction on the real algebraic number fields that Normaliz can use.

Normaliz has a library as well as a file interface. It can be reached from CoCoA, GAP \cite{GAP-NmzInterface}, Macaulay2, Singular, Python \cite{PyNormaliz} and SageMath. The full functionality is reached on Linux and Mac OS platforms, but the basic functionality for rational polyhedra is also available on MS Windows systems.

Its history goes back to the mid 90ies. For recent developments see \cite{BI2} and \cite{BSS}. The extension to algebraic polytopes was done in several steps since 2016. We are grateful to Matthias K\"oppe for suggesting it.

The work on algebraic polytopes has been done in cooperation with Vincent Delecroix (\eantic), Sebastian Gutsche (PyNormaliz), Matthias K\"oppe and Jean-Phiippe Labb\'e (integration into SageMath). A comprehensive article with these coauthors is in preparation.

\section{The icosahedron}
Let us specify the icosahedron, a Platonic solid, by its vertices:

\noindent\begin{minipage}[b]{0.5\textwidth}
\small
\begin{verbatim}
	
amb_space 3
number_field min_poly (a^2 - 5) embedding [2 +/- 1]
vertices 12
0 2 (a + 1) 4
0 -2 (a + 1) 4
2 (a + 1) 0 4
...
(-a - 1) 0 -2 4
Volume
LatticePoints
FVector
EuclideanAutomorphisms
	
\end{verbatim}
\end{minipage}
\hspace{1.5cm}
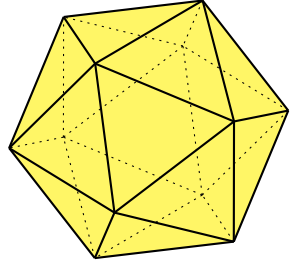
\begin{minipage}[t]{0.4\textwidth}
	\vspace*{-4.0cm}
	\begin{tikzpicture}%
	[x={(0.700041cm, -0.429565cm)},
	y={(0.714101cm, 0.419519cm)},
	z={(0.001418cm, 0.799673cm)},
	scale=2.00000,
	back/.style={dotted, thin},
	edge/.style={color=black!95!black, thick},
	facet/.style={fill=yellow,fill opacity=0.600000},
	vertex/.style={inner sep=0pt,circle,draw=black!25!black,fill=black!75!black,thick}]
	%
	%
	\coordinate (0.80902, 0.00000, 0.50000) at (0.80902, 0.00000, 0.50000);
	\coordinate (0.80902, 0.00000, -0.50000) at (0.80902, 0.00000, -0.50000);
	\coordinate (0.00000, 0.50000, 0.80902) at (0.00000, 0.50000, 0.80902);
	\coordinate (0.00000, 0.50000, -0.80902) at (0.00000, 0.50000, -0.80902);
	\coordinate (0.50000, 0.80902, 0.00000) at (0.50000, 0.80902, 0.00000);
	\coordinate (-0.50000, 0.80902, 0.00000) at (-0.50000, 0.80902, 0.00000);
	\coordinate (0.00000, -0.50000, 0.80902) at (0.00000, -0.50000, 0.80902);
	\coordinate (0.00000, -0.50000, -0.80902) at (0.00000, -0.50000, -0.80902);
	\coordinate (0.50000, -0.80902, 0.00000) at (0.50000, -0.80902, 0.00000);
	\coordinate (-0.80902, 0.00000, 0.50000) at (-0.80902, 0.00000, 0.50000);
	\coordinate (-0.80902, 0.00000, -0.50000) at (-0.80902, 0.00000, -0.50000);
	\coordinate (-0.50000, -0.80902, 0.00000) at (-0.50000, -0.80902, 0.00000);
	\fill[facet] (0.00000, -0.50000, 0.80902) -- (0.80902, 0.00000, 0.50000) -- (0.00000, 0.50000, 0.80902) -- cycle {};
	\fill[facet] (-0.80902, 0.00000, 0.50000) -- (0.00000, 0.50000, 0.80902) -- (0.00000, -0.50000, 0.80902) -- cycle {};
	\fill[facet] (0.50000, 0.80902, 0.00000) -- (0.80902, 0.00000, 0.50000) -- (0.00000, 0.50000, 0.80902) -- cycle {};
	\fill[facet] (0.50000, 0.80902, 0.00000) -- (0.80902, 0.00000, 0.50000) -- (0.80902, 0.00000, -0.50000) -- cycle {};
	\fill[facet] (0.50000, -0.80902, 0.00000) -- (0.80902, 0.00000, 0.50000) -- (0.00000, -0.50000, 0.80902) -- cycle {};
	\fill[facet] (0.50000, -0.80902, 0.00000) -- (0.80902, 0.00000, -0.50000) -- (0.00000, -0.50000, -0.80902) -- cycle {};
	\fill[facet] (0.50000, -0.80902, 0.00000) -- (0.80902, 0.00000, 0.50000) -- (0.80902, 0.00000, -0.50000) -- cycle {};
	\fill[facet] (-0.50000, -0.80902, 0.00000) -- (0.00000, -0.50000, 0.80902) -- (-0.80902, 0.00000, 0.50000) -- cycle {};
	\fill[facet] (-0.50000, -0.80902, 0.00000) -- (0.00000, -0.50000, 0.80902) -- (0.50000, -0.80902, 0.00000) -- cycle {};
	\fill[facet] (-0.50000, -0.80902, 0.00000) -- (0.00000, -0.50000, -0.80902) -- (0.50000, -0.80902, 0.00000) -- cycle {};
	\draw[edge,back] (0.80902, 0.00000, -0.50000) -- (0.00000, 0.50000, -0.80902);
	\draw[edge,back] (0.00000, 0.50000, 0.80902) -- (-0.50000, 0.80902, 0.00000);
	\draw[edge,back] (0.00000, 0.50000, -0.80902) -- (0.50000, 0.80902, 0.00000);
	\draw[edge,back] (0.00000, 0.50000, -0.80902) -- (-0.50000, 0.80902, 0.00000);
	\draw[edge,back] (0.00000, 0.50000, -0.80902) -- (0.00000, -0.50000, -0.80902);
	\draw[edge,back] (0.00000, 0.50000, -0.80902) -- (-0.80902, 0.00000, -0.50000);
	\draw[edge,back] (0.50000, 0.80902, 0.00000) -- (-0.50000, 0.80902, 0.00000);
	\draw[edge,back] (-0.50000, 0.80902, 0.00000) -- (-0.80902, 0.00000, 0.50000);
	\draw[edge,back] (-0.50000, 0.80902, 0.00000) -- (-0.80902, 0.00000, -0.50000);
	\draw[edge,back] (0.00000, -0.50000, -0.80902) -- (-0.80902, 0.00000, -0.50000);
	\draw[edge,back] (-0.80902, 0.00000, 0.50000) -- (-0.80902, 0.00000, -0.50000);
	\draw[edge,back] (-0.80902, 0.00000, -0.50000) -- (-0.50000, -0.80902, 0.00000);
	\draw[edge] (0.80902, 0.00000, 0.50000) -- (0.80902, 0.00000, -0.50000);
	\draw[edge] (0.80902, 0.00000, 0.50000) -- (0.00000, 0.50000, 0.80902);
	\draw[edge] (0.80902, 0.00000, 0.50000) -- (0.50000, 0.80902, 0.00000);
	\draw[edge] (0.80902, 0.00000, 0.50000) -- (0.00000, -0.50000, 0.80902);
	\draw[edge] (0.80902, 0.00000, 0.50000) -- (0.50000, -0.80902, 0.00000);
	\draw[edge] (0.80902, 0.00000, -0.50000) -- (0.50000, 0.80902, 0.00000);
	\draw[edge] (0.80902, 0.00000, -0.50000) -- (0.00000, -0.50000, -0.80902);
	\draw[edge] (0.80902, 0.00000, -0.50000) -- (0.50000, -0.80902, 0.00000);
	\draw[edge] (0.00000, 0.50000, 0.80902) -- (0.50000, 0.80902, 0.00000);
	\draw[edge] (0.00000, 0.50000, 0.80902) -- (0.00000, -0.50000, 0.80902);
	\draw[edge] (0.00000, 0.50000, 0.80902) -- (-0.80902, 0.00000, 0.50000);
	\draw[edge] (0.00000, -0.50000, 0.80902) -- (0.50000, -0.80902, 0.00000);
	\draw[edge] (0.00000, -0.50000, 0.80902) -- (-0.80902, 0.00000, 0.50000);
	\draw[edge] (0.00000, -0.50000, 0.80902) -- (-0.50000, -0.80902, 0.00000);
	\draw[edge] (0.00000, -0.50000, -0.80902) -- (0.50000, -0.80902, 0.00000);
	\draw[edge] (0.00000, -0.50000, -0.80902) -- (-0.50000, -0.80902, 0.00000);
	\draw[edge] (0.50000, -0.80902, 0.00000) -- (-0.50000, -0.80902, 0.00000);
	\draw[edge] (-0.80902, 0.00000, 0.50000) -- (-0.50000, -0.80902, 0.00000);
	\end{tikzpicture}
\end{minipage}
The first line specifies the dimension of the affine space. The second defines the unique positive sqare root of $5$ as the generator of the number field. It is followed by the $12$ vertices. Each of them is given as a vector with $4$ components for which the fourth component acts as a common denominator of the first three. Expressions involving $a$ are enclosed in round brackets. The last lines list the computation goals for Normaliz. (Picture by J.-P. Labb\'e)

Normaliz has a wide variety of input data types. For example, it would be equally possible to define the icosahedron by inequalities. Now we have a look into the output file. (We indicate omitted lines by \dots)
{\small
\begin{verbatim}
Real embedded number field:
min_poly (a^2 - 5) embedding [2.23606797...835961152572 +/- 5.14e-54]

1 lattice points in polytope
12 vertices of polyhedron
0 extreme rays of recession cone
20 support hyperplanes of polyhedron (homogenized)
f-vector:
1 12 30 20 1
embedding dimension = 4
affine dimension of the polyhedron = 3 (maximal)
rank of recession cone = 0 (polyhedron is polytope)
... 
volume (lattice normalized) = (5/2*a+15/2 ~ 13.090170)
volume (Euclidean) = 2.18169499062
Euclidean automorphism group has order 120
***********************************************************************
1 lattice points in polytope:
0 0 0 1
12 vertices of polyhedron:
...
0 extreme rays of recession cone:
20 support hyperplanes of polyhedron (homogenized):
(-a+1 ~ -1.236068) (-2*a+4 ~ -0.472136)                    0 1
...
  (a-1 ~ 1.236068)   (2*a-4 ~ 0.472136)                    0 1
\end{verbatim}
}
The output (in homogenized coordinates) is self-explanatory. Note that non-integral numbers in the output are printed as polynomials in $a$ together with a rational approximation. At the top we can see to what precision $\sqrt 5$ had to be computed. The automorphism group is described in another output file:
{
\small
\begin{verbatim}
Euclidean automorphism group of order 120
************************************************************************
3 permutations of 12 vertices of polyhedron
Perm 1: 1 2 4 3 7 8 5 6 10 9 11 12
Perm 2: 1 3 2 5 4 6 7 9 8 11 10 12
Perm 3: 2 1 3 4 6 5 8 7 9 10 12 11
Cycle decompositions 
Perm 1: (3 4) (5 7) (6 8) (9 10) --
Perm 2: (2 3) (4 5) (8 9) (10 11) --
Perm 3: (1 2) (5 6) (7 8) (11 12) --
1 orbits of vertices of polyhedron
Orbit 1 , length 12:  1 2 3 4 5 6 7 8 9 10 11 12
************************************************************************
3 permutations of 20 support hyperplanes
Perm 1: 2 1 5 6 3 4 7 8 11 12 9 10 13 14 17 18 15 16 20 19
...
Cycle decompositions 
Perm 1: (1 2) (3 5) (4 6) (9 11) (10 12) (15 17) (16 18) (19 20) --
...
1 orbits of support hyperplanes
Orbit 1 , length 20:  1 2 3 4 5 6 7 8 9 10 11 12 13 14 15 16 17 18 19 20
\end{verbatim}
}

\section{Computation goals for algebraic polyhedra}

The basic computation in linear convex geometry is the dualization of cones. We start from a cone $C\subset \RR^d$, given by generators $x_1,\dots,x_n$. The first (easy) step is to find a coordinate transformation that replaces $\RR^d$ by the vector subspace generated by $x_1,\dots,x_n$. In other words, we can assume $\dim C=d$.

The goal is to find a minimal generating set $\sigma_1,\dots,\sigma_s\in(\RR^d)^*$ of the dual cone $C^*=\bigl\{\lambda: \lambda(x_i)\ge 0, \ i=1,\dots,n\bigr\}$. Because of $\dim C=d$, the linear forms $\sigma_1,\dots,\sigma_s$ are uniquely determined up to positive scalars: they are the extreme rays of $C^*$. By a slight abuse of terminology we call the hyperplanes $S_i=\{x: \sigma_i(x)=0\}$ the \emph{support hyperplanes} of $C$.

Let $C_k$ be the cone generated by $x_1,\dots,x_k$. Normaliz proceeds as follows: 
\begin{enumerate}
\item It finds a basis of $\RR^d$ among the generators $x_1,\dots,x_n$, say $x_1,\dots,x_d$. Computing $C_d^*$ amounts to a matrix inversion. 
\item Iteratively it extends the cone $C_k$ to $C_{k+1}$, and shrinks $C_k^*$ to $C_{k+1}^*$, $k=d\dots,n-1$.
\end{enumerate}
Step 2 is done by \emph{Fourier-Motzkin elimination}: if $\sigma_1,\dots,\sigma_t$ generate $C_k^*$, then $C_{k+1}^*$ is generated by
$$
\bigl\{\sigma_i: \sigma_i(x_{k+1})\ge0 \bigr\}\cup \bigl\{\sigma_i(x_{k+1})\sigma_j-\sigma_j(x_{k+1})\sigma_i: \sigma_i(x_{k+1})\ >0, \sigma_j(x_{k+1}) < 0 \bigr\}.
$$
From this generating set of $C^*_{k+1}$ the extreme rays of $C^*_{k+1}$ must be selected.

This step is of critical complexity.  Normaliz has a sophisticated implementation in which \emph{pyramid decomposition} is a crucial tool; see \cite{BIS}. It competes very well with dedicated packages (see \cite{Koeppe}). The implementation is independent of the field of coefficients. As said above, $\RR$ can be replaced by an algebraic number  field $\AA$. In this case Normaliz uses the arithmetic over the field $\AA$ realized by \eantic, whereas arithmetic over $\QQ$ is avoided in favor of arithmetic over $\ZZ$.

In addition to the critical complexity caused by the combinatorics of cones, one must tame the coordinates of the linear combination $\lambda=\sigma_i(x_{k+1})\sigma_j-\sigma_j(x_{k+1})\sigma_i$. For example, if, over $\ZZ$, both $\sigma_i$ and $\sigma_j$ are divisible by $2$, then $\lambda$ is divisible by $4$. If this observation is ignored, a doubly exponential explosion of coefficients will happen. One therefore extracts the gcd of the coordinates. But there is usually no well-defined gcd of algebraic integers, and even if one has unique decomposition into prime elements, there is in general no Euclidean algorithm. Normaliz therefore applies two steps:
\begin{enumerate}
\item $\lambda$ is divided by the absolute value of the last nonzero component (or by another ``norm'').
\item All integral denominators are cleared by multiplication with their lcm.
\end{enumerate}
Computational experience has shown that these two steps together are a very good choice.

Normaliz tries to measure the complexity of the arithmetic in $\AA$ and to control the algorithmic alternatives of the dualization by the measurements. There are several ``screws'' that can be turned, and it is difficult to find the optimal tuning beforehand.

Normaliz computes lexicographic triangulations of algebraic cones in the same way as triangulations of rational cones. Their construction is interleaved with the extension from $C_k$ to $C_{k+1}$: the already computed triangulation of $C_k$ is extended by the simplicial cones generated by $x_{k+1}$ and those subcones in the triangulation of $C_k$ that are ``visible'' from $x_{k+1}$.

An algebraic polytope $P$ contains only finitely many integral points. They are computed by Normaliz' project-and-lift algorithm. The truncated Hilbert basis approaches, which Normaliz can also use for rational polytopes, are not applicable in the algebraic case. Once the lattice points are known, one can compute their convex hull, called the \emph{integer hull} of $P$.

At present Normaliz computes volumes only for full-dimensional algebraic polytopes.  The volume is the sum of the volumes of the simplices in a triangulation, and these are simply (absolute values of) determinants. We do not see any reasonable definition of ``algebraic volume'' for lower dimensional polytopes that could replace the lattice normalized volume. The latter is defined for all rational polytopes and is a rational number that can be computed precisely. 

It would certainly be possible to extend the computation of the approximate Euclidean volume to all algebraic polytopes, and this extension may be included in future Normaliz versions. Note that the Euclidean volume does in general not belong to $\AA$ if $P$ is lower dimensional. Its precise computation would require an extension of $\AA$ by square roots.

The computation of automorphism groups follows the suggestions in \cite{Bremner}. First one transforms the defining data into a graph, and then computes the  automorphism group of this graph by \emph{nauty} \cite{nauty}. For algebraic polytopes the Euclidean and the algebraic automorphism groups can be computed, and the combinatorial automorphism group is accessible for all polyhedra. 

The Euclidean automorphism group is the group of rigid motions of the ambient space that map the polytope to itself, and the algebraic automorphism group is the group of affine transformations over $\AA$ stabilizing the polytope. Both groups are finite, as well as the combinatorial automorphism group, the automorphism group of the face lattice, which can be computed from the facet-vertex incidence vectors, just as in the rational case.

We do not try to define the algebraic (or Euclidean) automorphism group for unbounded polyhedra. First of all, the algebraic automorphism group is infinite in general. Second, it would have to be realized as the permutation group of a vector configuration, and there seems to be no reasonable way to norm the involved vectors. But for polytopes we can and must use the vertices.
 
\section{Scaled convex hull computations}

We illustrate the influence of the algebraic number field on the computation time by some examples. For each of them we start from a cone (over a polyhedron) that is originally defined over the integers. Then we scale some coordinates by elements of the field $\AA$. This transformation preserves the combinatorial structure throughout. It helps to isolate the complexity of the arithmetic operations. The types of arithmetic that we compare are 
\begin{quote}
	int: original input, computation with machine integers,\\
	mpz: same input as int, but computation with GMP mpz\_class integers,\\
	rat: same input as int, but computation in $\QQ[\sqrt 5]$,\\
	sc2: scaled input in $\QQ[\sqrt 5]$,\\
	sc8: scaled input in $\QQ[\sqrt[8] 5]$,\\
	p12: scaled input in $\QQ[a]$, $a^{12} + a^6+a^5+a^2- 5=0$, $a>1$.
\end{quote}

The test candidates are A553 (from the Ohsugi-Hibi classification of contingency tables \cite{OH}), the cone q27f1 from \cite{Koeppe}, the linear order polytope for $S_6$, and the cyclic polytope of dimension $15$ with $30$ vertices. The last two are classical polytopes. While the other three cones are given by their extreme rays, q27f1 is defined by $406$ equations and inequalities.

\begin{table}[hbt]
\caption{Combinatorial data of the test candidates}\label{comb}
\centering
\small
\tabcolsep1ex
\begin{tabular}{| r | r  | r | r | r |}
\hline
\strut & amb\_space & dim& ext rays & supp hyps\\
\hline
\strut A553& $55$ & $43$ & $75$ & $306,955$\\
\hline
\strut q27f1& $30$ & $13$ & $68,216$ & $92$ \\
\hline
\strut lo6& $16$& $16$&  $720$ & $910$ \\
\hline
\strut cyc15-30 &$16$ & $16$ & $30$ & $341088$\\
\hline
\end{tabular}
\end{table}

The Normaliz version is 3.8.4, compiled into a static binary with gcc 5.4 under Ubuntu 16-04. The computations use $8$ parallel threads (the default choice of Normaliz). They were taken on the author's PC with an AMD Ryzen 7 1700X at 3.2 GHz. Table \ref{perf} lists wall times in seconds. As a rule of thumb, for a single thread the times must be multiplied by $6$.
 
\begin{table}[hbt]
\caption{Wall times of scaled convex hull computations in seconds}\label{perf}
\centering
\small
\tabcolsep1ex
\begin{tabular}{|l|r|r|r|r|r|}
\hline
\strut coeff    & A553       & q27f1  & lo6   & cyc15-30 \\
\hline
\strut int        &57        &16      & 5     & --\\    
\hline                       
\strut mpz        &299       &58      & 5     & 7\\    
\hline                       
\strut rat        &277       &40      & 5     & 7 \\  
\hline                       
\strut sc2        &783       &166     & 4     & 14 \\ 
\hline                       
\strut sc8        &1272      &475     &15     & 28\\
\hline                       
\strut p12        &2908      &905     &31     & 42\\                    
\hline
\end{tabular}
\end{table}

The cyclic polytope and all intermediate polytopes coming up in its computation are simplicial. Therefore it profits from Normaliz' special treatment of simplicial facets---almost everything can be done by set theoretic operations. Also lo6 is combinatorially not complicated.  That lo6 is fastest with sc2, is caused by the fine tuning of the pyramid decomposition, which is not always optimal.

Surprisingly, rat is faster than mpz for A553 and q27f1. This can be explained by the fact that linear algebra over $\ZZ$ must use the Euclidean algorithm, and therefore needs more steps than the true rational arithmetic of rat.

\end{document}